\documentclass[14pt,a4paper]{article}

\usepackage[utf8]{inputenc}
\usepackage[english,russian]{babel}
\usepackage{indentfirst}
\usepackage{misccorr}
\usepackage{graphicx}
\usepackage{amsmath}
\usepackage{amssymb}
\usepackage{amsthm}
\usepackage{hyperref}
\usepackage{color}

\newif\ifhide
\hidefalse

\newtheorem{theorem}{Теорема}
\newtheorem{Lemma}[theorem] {Лемма}
\newtheorem{Proposition}[theorem]{Предложение}

\newtheorem{Corollary}[theorem] {Следствие}

\title{~\hspace{-77mm}{\small УДК: 519.21:519.65:517.518.8}
\\~\\
О частных производных модифицированных полиномов Бернштейна - Станку для функций нескольких переменных
\footnote{Для обоих авторов 
данная работа была поддержана грантом Фонда развития теоретической физики и математики «БАЗИС».}}

\author{А.Ю.Веретенников\footnote{Институт проблем передачи информации им. А.А.Харкевича, email: ayv@iitp.ru}, Н.М.Мазутский\footnote{МГУ им. М.В.Ломоносова, 
 email: 
 nikolaymazutsky@gmail.com}}

\begin{document}

\maketitle
\begin{abstract}
\noindent
Целью работы является доказательство аппроксимации смешанных производных второго порядка для функции нескольких переменных  в норме $L_1$ такими же производными модифицированных  полиномов Бернштейна при минимальной возможной гладкости. 

\medskip

\noindent
Ключевые слова: модифицированные полиномы Бернштейна; аппроксимации смешанных производных

\medskip


\end{abstract}

\section{Введение}
Для нового (на тот момент) доказательства аппроксимационной полиномиальной теоремы Вейерштрасса на основе закона больших чисел (ЗБЧ) в схеме Бернулли С.Н.Бернштейн \cite{Bernstein} предложил полиномы, названные затем его именем:
$$
{B}_{n}(f,x)= \sum\limits_{k=0}^{n} f\left( \frac{k}{n}\right) \dbinom{n}{k}  x^{k} {(1-x)}^{n-k}, \quad 0\le x\le 1.
$$
Кажущаяся на первый взгляд искусственной интерпретация здесь переменной $x$ как вероятности и применение ЗБЧ приводит к теореме о сходимости для любой функции $f\in C[0,1]$:
$$
\sup_{\vert x \vert \le 1}\vert B_n(f,x) - f(x)\vert \to 0, \quad n\to\infty.
$$
Дальнейшее развитие этой идеи шло по разным направлениям: как оценить скорость сходимости при той, или иной гладкости приближаемой функции; можно ли производными полиномов $B_n$ приблизить соответствующие производные функции $f$; что можно сказать в отсутствие предположения о непрерывности $f$; наконец, существуют ли аналоги полиномов $B_n(x)$ для приближения функций нескольких переменных на кубе, или на симплексе, и как обстояит дело с приближением частных производных $f$ в этом случае. Различными авторами были предложены разнообразные модификации полиномов Бернштейна, см., в частности, \cite{Lorentz 2, Stancu, Voronovskaya, Vedinskiy, Pendina, Polovinkina, Kirov, Kwun, Kantorovich}.  В работе \cite{Veretennikov} было установлено, что естественные аналоги полиномов Бернштейна приближают функции нескольких переменных на кубе и  на симплексе вместе с их производными определенного порядка, или индекса, в предположении о непрерывности производных этого порядка, или индекса (то есть, не обязательно {\it всех} производных, скажем, второго порядка, а лишь тех, которые определены и непрерывны). 
В настоящей работе на основе некоторой модификации полиномов Бернштейна - Станку (см. \cite{Stancu}) исследован многомерный (формально двумерный) случай и для него -- вопрос о приближении смешанной частной производной второго порядка в интегральной норме и без предположения о ее непрерывности. Для полноты изложения сформулирован и доказан также предварительный результат о приближении в той же норме частных производных первого порядка. Работа состоит из шести разделов: Введение, О более ранних результатах, Основной результат (теорема \ref{thm1}), Вспомогательный результат, Доказательство теоремы \ref{thm1} и Замечания о (равенстве) смешанных производных. Все результаты, леммы и т.п. имеют сквозную нумерацию ради удобства их быстрого нахождения при чтении. 

\section{
О более ранних результатах}
Перед  презентацией основного результата статьи напомним  основополагающую теорему Бернштейна \cite{Bernstein}, а также результаты Хлодовского  \cite{Khlodovskii} и Канторовича \cite{Kantorovich}, а также еще одно, не столь давнее, предложение \ref{Veretennikov}. Все они, кроме последнего,  относятся к приближению функций одной переменной на $[0,1]$, тогда как теорема \ref{thm1} (см. ниже) -- к приближению функции (и ее производных) двух переменных; однако, 
это может помочь читателю быстро войти в курс дела.  Все ``чужие'' (то есть, взятые из других статей) теоремы будем называть предложениями. Поскольку работы \cite{Khlodovskii} и \cite{Kantorovich} нельзя назвать легко доступными, см. формулировки их результатов в  \cite[Теорема 1]{Lorentz} и \cite[Теорема 2.1.1]{Lorentz 2}, соответственно. 

\begin{Proposition}[Бернштейн {\cite{Bernstein}}]
Для любой функции $f$, непрерывной на $[0,1]$, последовательность полиномов $B_n(f,x)$ при $n \to \infty$ сходится к $f(x)$ равномерно на $[0,1]$.\\
\end{Proposition}

\begin{Proposition}[Хлодовский \cite{Khlodovskii}, Лоренц {\cite[Теорема 1]{Lorentz}}]\label{khlod}
Если для ограниченной функции $f(x)$ в точке $x_0$ существует r-тая производная $f^{(r)}(x_0)$, то 
$$
\displaystyle B_n^{(r)}(f,x_0) \xrightarrow{} f^{(r)}(x_0), \quad n \to \infty.
$$
\end{Proposition}
Ссылка \cite{Khlodovskii} дает лишь название устной презентации автора на Всесоюзном съезде математиков; тексты большей части докладов не публиковались. Тем не менее, в результате этого выступления  данный результат (предложение \ref{khlod}) затем получил известность в литературе под именем своего автора.

В работе \cite{Kantorovich} (см. также изложение в \cite[глава 2]{Lorentz 2})  в качестве видоизменения полиномов для исследования аппроксимационных свойств нерегулярных функций было предложено заменить интегрируемую на $[0,1]$ функцию $\displaystyle f(x)$ на $\displaystyle F(x)=\int_{0}^{x}f(t)dt.$ Таким образом, {\it производная} полинома для функции $F(x)$ примет следующий вид:\\
$$
\displaystyle P_n(f,x):= \frac{d}{dx}B_{n+1}(F,x)=\sum\limits_{k=0}^{n}\dbinom{n}{k}x^{k} {(1-x)}^{n-k}(n+1)\int_{\frac{k}{n+1}}^{\frac{k+1}{n+1}}f(t)dt.
$$

\begin{Proposition}[Канторович \cite{Kantorovich}]
В любой точке $x \in (0,1)$, где $f(x)$ -- производная от интеграла $\displaystyle F(x)=\int_{0}^{x}f(t)dt$, то есть почти всюду, имеет место сходимость
$$\displaystyle P_n(f,x) \xrightarrow{} f(x), \quad n \to \infty$$
\end{Proposition}
\noindent
Изложение этого результата можно также найти в \cite[Теорема 2.1.1]{Lorentz 2}, со ссылкой на Канторовича.

Одномерные аналоги введенных далее в (\ref{Bt}) модифицированных полиномов $\tilde B_n$ идеологически были бы близки данной модификации Канторовича; однако, тип сходимости, утверждаемой в теореме \ref{thm1} ниже, несколько иной. Кроме того, вспомним, что в многомерном случае нет полного аналога одномерной теоремы Лебега о дифференцировании, так что, на самом деле, автоматически  результат из \cite{Kantorovich} на многомерный случай не переносится. Также подчеркнем, что целью данной работы является иная задача, а именно -- приближение смешанной производной второго порядка.

Один из двух классических вариантов полиномов Бернштейна в случае функции двух переменных имеет следующий вид:
$$
\displaystyle  {B}_{n}(f,x_1,x_2)= \sum\limits_{k_1, k_2=0}^{n} f\left( \frac{k_1}{n},\frac{k_2}{n}\right) \dbinom{n}{k_1} \dbinom{n}{k_2} x_{1}^{k_1} {(1-x_1)}^{n-k_1} x_{2}^{k_2} {(1-x_2)}^{n-k_2}, 
$$
где $0\le x_1, x_2 \le 1$, и аналогично он может быть задан на $d$-мерном кубе $0\le x_1, \ldots, x_d \le 1$. Именно этот вид будет исследован далее в работе. (Второй известный вариант задает полином и приближает функцию не на квадрате/кубе, а на симплексе; в одномерном случае обе версии, естественно, совпадают.) 
В работе \cite{Veretennikov} был установлен (см. {\cite[Теорема 4]{Veretennikov}}) факт равномерной сходимости на $d$-мерном кубе $m$-й производной полинома Бернштейна многих переменных к $m$-й производной функции $f\in C^m(\mathbb{R}^d)$: $\displaystyle {B}_{n}^{(m)}
	\left( f,x\right) \overrightarrow{\xrightarrow{}} f^{(m)}\left( x\right), \quad n \to \infty$, а также аналогичное утверждение доказано и для мультииндекса $k=(k_1,..,k_d): \displaystyle {B}_{n}^{(k)}
	\left( f,x\right) \overrightarrow{\xrightarrow{}} f^{(k)}\left( x\right), \quad n \to \infty$, в предположении существования лишь  непрерывной производной, соотвествующей этому мультииндексу, но не обязательно всех производных данного порядка. То, что данный факт может представлять самостоятельный интерес, демонстрируется теоремой Г.П.Толстова (см. ниже теорему \ref{thm2}), где утверждается не сходимость, а лишь равенство смешанных производных взятых в различных порядках, однако в условия входит требование существования {\bf всех} производных порядка 2; такое предположение характерно и для  других версий этого результата; ни в предложении \ref{Veretennikov}, ни в теореме \ref{thm1} ниже оно не накладывается.

\begin{Proposition}[{\cite[Теорема 4]{Veretennikov}}]\label{Veretennikov}
	1) Если $\displaystyle m>0, f \in C^{m}({\mathbb{R}}^d),$
	то\\ \\ $\displaystyle B_{n}^{(m)}
	\left( f,x\right) \overrightarrow{\xrightarrow{}} f^{(m)}\left( x\right), n \to \infty,
	x \in {{K}}^{d},$ где ${{K}}^{d} - $ d-мерный куб.\\ \\
	2) Если ${\mathbf k}=(k_1,...,k_d)$ -- мультииндекс и $f \in C^{{\mathbf k}}({\mathbb{R}}^d)$,
	то\\ \\ $\displaystyle B_{n}^{({\mathbf k})}
	\left( f,x\right) \overrightarrow{\xrightarrow{}} f^{({\mathbf k})}\left( x\right), n \to \infty,
	x \in {{K}}^{d},$ где ${{K}}^{d} - $ d-мерный куб.
\end{Proposition}

По сравнению с последним указанным результатом, в теореме \ref{thm1} ниже  удалось отказаться от условия непрерывности смешанной производной второго порядка функции $f$ на $[0,1]^2:=K$, и при условии $\displaystyle\frac{\partial^2 f}{\partial x_1 \partial x_2}\in L_1(K)$ и некоторых других установлен факт сходимости смешанной частной производной второго порядка по переменным $x_1, x_2$ модифицированного полинома Бернштейна -- Станку к смешанной частной производной второго порядка функции $f$ в норме $L_1$ по $(x_1,x_2)\times\omega$ относительно прямого произведения двумерной меры Лебега на $K$ и вероятностной  меры $\mathsf P$ на измеримом пространстве $(\Omega, {\cal F}, \mathsf P)$ (см. точные определения в следующем разделе). Здесь $(x_1,x_2) \in K$, $\omega \in \Omega$. Отметим, что это иной тип сходимости и по сравнению с вариантом сходимости в указанном выше результате из \cite{Veretennikov}. Предложение \ref{Veretennikov} будет использована далее в ходе доказательства теоремы \ref{thm1} в разделе 5;  по этой причине и ради удобства читателя  этот результат и сформулирован выше ради удоства ссылок. Далее в работе рассматриваются лишь функции двух переменных, поскольку предметом изучения являются смешанные производные второго порядка, для которых большее число переменных не актуально.

\section{Основной результат}
Для формулировки основного результата введем две независимые случайные величины\footnote{Напомним, что $U[0,1]$ является стандартным обозначением равномерно распределенных на $[0,1]$ случайных величин; знак $\mathsf E$ означает математическое ожидание по мере $\mathsf P$.} $\alpha_i \in U[0,1]$, $i = 1,2$, на вероятностном пространстве $(\Omega, {\cal F}, \mathsf P)$, равномерно распределенные на $[0,1]$, и положим $\alpha = (\alpha_1, \alpha_2)$. 
Полиномы {\em типа Бернштейна -- Станку,} исследуемые в настоящей работе, будут иметь следующий вид:
\begin{align}\label{Bt}
\!\!\displaystyle \tilde {B}_{\alpha,n}(f,x_1,x_2)\! = &  \!\!\!\!\!\sum\limits_{k_1, k_2=0}^{n} \!\!f\left( \frac{k_1+\alpha_1}{n},\frac{k_2+\alpha_2}{n}\right)\!\!\! \dbinom{n}{k_1}\! \dbinom{n}{k_2} 
  \nonumber \\\\ \nonumber
&\times x_{1}^{k_1} {(1\!-\!x_1)}^{n-k_1} x_{2}^{k_2} {(1\!-\!x_2)}^{n-k_2}.
\end{align}
Полезно доопределить $\displaystyle f(x_1,x_2)=0$ вне $K=[0,1]^2$, что и сделаем. 
Мы определили полиномы на квадрате, так как нас будет интересовать именно двумерный вариант этих полиномов; аналогичным образом они могут быть определены и на $d$-мерном кубе. 
Как уже отмечалось выше, похожие полиномы с неслучайными $\alpha_i$ были впервые  предложены Д.Д.Станку в \cite{Stancu}.

\begin{theorem}\label{thm1}    
Пусть $f(x_1,x_2)$ --- ограниченная, борелевская на $[0,1]^2:=K$.
Пусть существуют (классические) производные $\displaystyle\frac{\partial{f(x_1,x_2)}}{\partial{x_1}}\in L_1(K) \; \text{и}$ $\displaystyle \frac{\partial{f(x_1,x_2)}}{\partial{x_2}} \in L_1(K)$, и также что существует смешанная производная второго порядка 
$f_{x_1 x_2}(x_1, x_2)= \displaystyle\frac{\partial^2 f(x_1,x_2)}{\partial x_1 \partial x_2} \equiv \frac{\partial}{\partial x_1} \left(\frac{\partial f(x_1,x_2)}{ \partial x_2}\right)\in L_1(K)$. 
Тогда 
\begin{equation}\label{T1-1}
\displaystyle {\mathsf E}_{} \int\limits_0^1 dx_1 \int\limits_0^1 dx_2 \bigg| \frac{\partial^2 \widetilde{B}_{\alpha,n}(f,x_1,x_2)}{\partial x_1 \partial x_2} - \frac{\partial^2 f(x_1,x_2)}{\partial x_1 \partial x_2}\bigg| \xrightarrow{} 0, \quad n \to \infty.
\end{equation}
\end{theorem}

\begin{Corollary}
Если в дополнение к условиям части 2 теоремы \ref{thm1} предположить, что существует также частная смешанная производная $f_{x_2x_1}\in L_1(K)$, то аналогичная (\ref{T1-1}) сходимость имеет место и для нее и имеем место равенство п.в. $(x_1,x_2)$:
\[
f_{x_2x_1}(x_1,x_2) = f_{x_1x_2}(x_1,x_2).
\]
\end{Corollary}

\section{Вспомогательные результаты}
По мнению авторов, следующая лемма хорошо известна в теории интегрирования.  Однако, найти точную ссылку нам не удалось, и потому, по настоянию рецензента, приводим ее с доказательством, естественно, не претендуя на авторство\footnote{С благодарностью примем совет, содержащий нужную ссылку, и тогда доказательство можно будет убрать.}. 

\begin{Lemma}\label{lemma0}
{1.} Пусть $g\in L_1([0,1])$; вне отрезка $[0,1]$ доопределим функцию $g$ нулем. Тогда 
$$
\int_0^1 |g(x+t) - g(x)|dx \to 0, \quad t\to 0.
$$

\noindent
{2.}
Пусть $g\in L_1([0,1]^2)$; вне квадрата $K=[0,1]^2$ доопределим функцию $g$ нулем. Тогда 
$$
\iint\limits_K |g(x_1+t_1, x_2+t_2) - g(x_1,x_2)|dx_1dx_2 \to 0, \quad t_1,t_2\to 0.
$$
\end{Lemma}

\begin{proof}
1. Пусть $\delta>0$, и выберем любую функцию $h\in C([0,1])$ такую, что 
$$
\|g-h\|_{L_1([0,1])} \le \delta.
$$
Тогда оцениваем, 
\begin{align*}
&\int_0^1 |g(x+t) - g(x)|dx \le \int_0^1 |g(x+t) - h(x+t)|dx 
 \\\\
&+ \int_0^1 |g(x) - h(x)|dx + \int_0^1 |h(x+t) - h(x)|dx 
 \\\\
& \le 2\delta +  \int_0^1 |h(x+t) - h(x)|dx.
\end{align*}
Здесь последний интеграл стремится к нулю при $t\to 0$ в силу поточечной сходимости $h(x+t)$ к $h(x)$ и посольку функция $h$ ограничена, по теореме Лебега об ограниченной сходимости. 

\medskip

\noindent
2. Пусть вновь $\delta>0$, и выберем любую функцию $h\in C(K)$ такую, что 
$$
\|g-h\|_{L_1(K)} \le \delta.
$$
Подчеркнем, что на этот раз обе функции $g$ и $h$ заданы на $K$. 
Тогда оцениваем, 
\begin{align*}
&\iint\limits_K |g(x_1+t_1, x_2+t_2) - g(x_1,x_2)|dx_1dx_2
 \\\\ & 
\le  \iint\limits_K |g(x_1+t_1, x_2+t_2) - h(x_1+t_1,x_2+t_2)|dx_1dx_2
 \\\\
&+ \iint\limits_K |g(x_1, x_2) - h(x_1,x_2)|dx_1dx_2 
+\iint\limits_K |h(x_1+t_1, x_2+t_2) - h(x_1,x_2)|dx_1dx_2
 \\\\
&  \le 2\delta +   \iint\limits_K |h(x_1+t_1, x_2+t_2) - h(x_1,x_2)|dx_1dx_2.
\end{align*}
Как и в одномерном случае, в силу (даже равномерной) непрерывности и ограниченности $h$ на $K$, последний интеграл стремится к нулю при $t_1,t_2\to 0$, по теореме Лебега об ограниченной сходимости. Лемма \ref{lemma0} доказана. 
\end{proof}

\medskip

Так как в дальнейшем пойдет речь о производных функций, уже заданных достаточно громоздкими выражениями, а формулы для производных будут, естественно, еще более громоздки, то полезно ввести следующие обозначения, существенно сокращающие размеры выкладок, хоть они все равно  остаются не маленькими:
\\ \\ $\Delta_{(x_1)}f(x_1,x_2)=\Delta_{\frac{1}{n},x_1}f(x_1,x_2):=f(x_1+\frac{1}{n}, x_{2})-f(x_1, x_2).$\\
\newline
$\Delta_{(x_2)}f(x_1,x_2)=\Delta_{\frac{1}{n},x_2}f(x_1,x_2):=f(x_1, x_{2}+\frac{1}{n})-f(x_1, x_2).$\\
\newline
$\displaystyle f_{x_i}:=\frac{\partial{f(x_1,x_2)}}{\partial{x_i}}, \quad f_{x_i,x_j}:=\frac{\partial^2{f(x_1,x_2)}}{\partial{x_i}\partial{x_j}} = 
\frac{\partial{}}{\partial{x_i}} \left( \frac{\partial{f(x_1,x_2)}}{\partial{x_j}} \right), \quad i,j=1, 2.$

\begin{Lemma}\label{lemma1}
1. Пусть функция $f$ ограниченная и борелевская на $K$, и при всех $x_1, x_2$ существует классическая производная первого порядка $\displaystyle \quad \frac{\partial f(x_1, x_2
\newline)}{\partial x_1} \in L_1(K)$.
Тогда
\begin{equation}\label{eqle1}
\!\!\int_{0}^{1}\!dx_1 \!\int_{0}^{1}\!dx_2 \left| n \left( f_{}\left(x_1\!+\!\frac{1}{n}, x_2\right)
\!-\!f_{}\bigg( x_1, x_2\bigg) \right) \!-\! \frac{\partial f(x_1, x_2)}{\partial x_1 }\right|
\xrightarrow{} 0, \; n \to \infty.
\end{equation}
Аналогичное утверждение справедливо и для производной по другой переменной $f_{x_2}$ при условии ее существования и интегрируемости $\displaystyle \quad \frac{\partial f(x_1, x_2
\newline)}{\partial x_2} \in L_1(K)$:
\begin{equation}\label{le6-2}
    \displaystyle \int_{0}^{1}dx_1 \int_{0}^{1}dx_2 \left| n \left( f_{}\left( x_1, x_2+\frac{1}{n}\right)
    -f_{}\bigg( x_1, x_2\bigg) \right) - \frac{\partial f(x_1, x_2)}{\partial x_2 }\right|
    \xrightarrow{} 0, \quad n \to \infty.
\end{equation}

\medskip

\noindent
2.	Пусть при всех $x_1, x_2$ определена классическая смешанная производная $\displaystyle \frac{\partial^2 f(x_1, x_2
\newline)}{\partial x_1 \partial x_2} \in L_1(K)$ (в частности, функция $f_{x_2}( x_1, x_2)$ ограничена при всяком $x_2$). Тогда   
\begin{equation}\label{eqle2}
 \int_{0}^{1}dx_1 \int_{0}^{1}dx_2 \left| n \left( f_{x_2}\left( x_1+\frac{1}{n}, x_2\right)
    -f_{x_2}\bigg( x_1, x_2\bigg) \right) - \frac{\partial^2 f(x_1, x_2)}{\partial x_1 \partial x_2}\right|
    \xrightarrow{} 0, \quad n \to \infty.
\end{equation}
Другими словами,
$$
\int_{0}^{1}dx_1 \int_{0}^{1}dx_2 \left\vert
n\Delta_{(x_1)}f_{x_2}(x_1,x_2) - f_{x_1 x_2}(x_1,x_2)\right\vert \to 0, \quad n \to \infty.
$$
Также 
\begin{equation}\label{dxxd}
\Delta_{(x_1)}f_{x_2}(x_1,x_2) = \frac{\partial }{\partial x_2}\Delta_{(x_1)}f_{}(x_1,x_2),
\end{equation} 
и 
\begin{equation}\label{eqle3}
\int_{0}^{1}dx_1 \int_{0}^{1}dx_2 \left\vert
n^2\Delta_{(x_2)}\Delta_{(x_1)}f_{}(x_1,x_2) \! -\!
n\Delta_{(x_1)}f_{x_2}(x_1,x_2)\right\vert \to 0, \; n \to \infty.
\end{equation}
Аналогично, если при всех $x_1, x_2$ определена классическая смешанная производная $\displaystyle \frac{\partial^2 f(x_1, x_2
\newline)}{\partial x_2 \partial x_1} \in L_1(K)$ (и, в частности, функция $f_{x_1}( x_1, x_2)$ ограничена при всяком $x_1$), то при 
$n \to \infty$ имеет место сходимость 
\begin{equation}\label{le6-5}
\int_{0}^{1}dx_1 \int_{0}^{1}dx_2 \left| n \left( f_{x_1}\left( x_1+\frac{1}{n}, x_2+\frac{1}{n}\right)
-f_{x_1}\bigg( x_1, x_2\bigg) \right) - \frac{\partial^2 f(x_1, x_2)}{\partial x_2 \partial x_1}\right|
\xrightarrow{} 0,
\end{equation}
и также
\begin{equation}\label{le6-6}
\!\int_{0}^{1}\!dx_1 \!\int_{0}^{1}\!dx_2 \left\vert
n^2\Delta_{(x_1)}\Delta_{(x_2)}f_{}(x_1,x_2)\! -\!
n\Delta_{(x_2)}f_{x_1}(x_1,x_2)\right\vert \to 0, \; n \to \infty.
\end{equation}
\end{Lemma}

\begin{proof}
В доказательстве неоднократно будут использованы вспомогательные функции $g\in C(K)$ для приближения функции $f$ и ее производных, а также $F^\varepsilon$ для некоторых интегральных разностей. Для различных $f, f_{x_1}, f_{x_2}, f_{x_1x_2}$ и для различных разностей эти функции будут, вообще говоря, всякий раз различными. 

\medskip

1. Установим справедливость утверждения (\ref{eqle1}).  С этой целью рассматрим приближения функции $f$ гладкими $f^\varepsilon \in C^\infty(\mathbb R^2)$. 
Будем для определенности считать, что $f^{\varepsilon}$ является  функцией, приближающей функцию $f$ посредством свертки с помощью одного и того же неотрицательного ядра $\varphi \in C_0^{\infty}(\mathbb R)$ (где $0$ означает компактность носителя) по каждой из двух переменных, где  $\displaystyle \int_{-\infty}^{\infty}\varphi(x)dx=1$, и $\displaystyle \varphi^{\varepsilon}(x)=\frac{1}{\varepsilon}\varphi(\frac{x}{\varepsilon})$. Тогда, как известно, 
$$
f^{\varepsilon}(x_1, x_2)=\int dy_1\int dy_2 \left[f(x_1-y_1, x_2-y_2) \varphi^{\varepsilon}(y_1)\varphi^{\varepsilon}(y_2)\right] \in C^{\infty}.
$$ 
В дальнейшем, при появлении произведения $\varphi^{\varepsilon}(y_1)\varphi^{\varepsilon}(y_2)$ будем иногда использовать для краткости выражения
$$
\varphi^{\varepsilon}(y_1)\varphi^{\varepsilon}(y_2) = : \bar \varphi^{\varepsilon}(y_1,y_2) \quad \text{и} \quad f^\varepsilon = f*\bar \varphi^\varepsilon.
$$

\medskip

Покажем, что в условиях данного пункта леммы оператор производной по $x_1$ перестановочен с интегрированием, то есть, что   
\begin{equation}\label{partfe}
\partial_{x_1}f^{\varepsilon}(x_1, x_2)\equiv
f_{x_1}^{\varepsilon}(x_1, x_2)=\int dy_1\int dy_2 \left[f_{x_1}(x_1-y_1, x_2-y_2) \varphi^{\varepsilon}(y_1)\varphi^{\varepsilon}(y_2)\right].
\end{equation}
Для этого достаточно установить, что при всех $0\le x_1 \le x'_1 \le 1$ и при любом $x_2$ справедливо равенство
$$
f^{\varepsilon}(x'_1, x_2) - f^{\varepsilon}(x_1, x_2) 
= \int_{x_1}^{x'_1} dz \int dy_1\int dy_2 \left[f_{z}(z-y_1, x_2-y_2) \varphi^{\varepsilon}(y_1)\varphi^{\varepsilon}(y_2)\right].
$$ 
Имеем, при условии, что все интегралы абсолютно сходятся (а так оно и есть в силу условий леммы),
\begin{align*}
&f^{\varepsilon}(x_1, x_2) 
+ \int_{x_1}^{x'_1} dz \int dy_1\int dy_2 \left[f_{z}(z-y_1, x_2-y_2) \varphi^{\varepsilon}(y_1)\varphi^{\varepsilon}(y_2)\right]
 \\\\
&= \int dy_1\int dy_2 \left[f(x_1-y_1, x_2-y_2) \varphi^{\varepsilon}(y_1)\varphi^{\varepsilon}(y_2)\right] 
 \\\\
&+  \int dy_1\int dy_2 \int_{x_1}^{x'_1} dz \left[f_{z}(z-y_1, x_2-y_2) \varphi^{\varepsilon}(y_1)\varphi^{\varepsilon}(y_2)\right] 
 \\\\
& = \iint   \left[f(x_1-y_1, x_2-y_2) +\int_{x_1}^{x'_1} dz f_{z}(z-y_1, x_2-y_2)  \right] 
\varphi^{\varepsilon}(y_1)\varphi^{\varepsilon}(y_2)dy_1 dy_2
 \\\\
& = \iint   f(x'_1-y_1, x_2-y_2) \varphi^{\varepsilon}(y_1)\varphi^{\varepsilon}(y_2)dy_1 dy_2 = f^{\varepsilon}(x'_1, x_2),
\end{align*}
по теореме Фубини и в силу формулы Ньютона - Лейбница, что и требовалось.

\medskip

Далее, оценим верхний предел при $n\to\infty$ (неотрицательного) интеграла в левой части (\ref{eqle1}), зависящего  от $n$,  но не от $\varepsilon$, c помощью вспомогательной функции $f^\varepsilon$. Покажем, что этот верхний предел не превосходит некоторого (также неотрицательного) выражения от $\varepsilon$, которое, в свою очередь, стремится к нулю при $\varepsilon \to 0$. Другими словами, при оценках вспомогательных выражений, зависящих от $n$ и $\varepsilon$, нас будут интересовать повторные пределы вида $\lim_{\varepsilon\to 0} \limsup_{n\to\infty}$. (Если бы выражения под интегралом не были неотрицательны, то, конечно, надо было бы исследовать также пределы вида $\lim_{\varepsilon\to 0} \liminf_{n\to\infty}$, однако, в данном случае это лишнее.)
Имеем,
\begin{align*}
&\int_{0}^{1}dx_1 \int_{0}^{1}dx_2 \bigg| n \left( f\left( x_1+\frac{1}{n}, x_2\right)
-f\bigg( x_1, x_2\bigg) \right) - f_{x_1}\bigg(x_1, x_2\bigg)\bigg|
 \\ \\
&\leqslant\int_{0}^{1}dx_1 \int_{0}^{1}dx_2\bigg|f^{\varepsilon}_{x_1}\bigg(x_1, x_2\bigg)-f_{x_1}\bigg(x_1, x_2\bigg)\bigg|
 \\ \\
&+\int_{0}^{1}dx_1 \int_{0}^{1}dx_2\bigg|n \left( f^{\varepsilon}\left( x_1+\frac{1}{n}, x_2\right)-f^{\varepsilon}\bigg( x_1, x_2\bigg) \right) - f^{\varepsilon}_{x_1}\bigg(x_1, x_2\bigg)\bigg|
 \\ \\
&+\int_{0}^{1}dx_1 \int_{0}^{1}dx_2\bigg|n \left( f\left( x_1+\frac{1}{n}, x_2\right)
-f\bigg( x_1, x_2\bigg) \right)
 \\ \\
&-n \left( f^{\varepsilon}\left( x_1+\frac{1}{n}, x_2\right)-f^{\varepsilon}\bigg( x_1, x_2\bigg) \right)\bigg| =:J_{1}^{\varepsilon}+J_{2}^{n,\varepsilon}+J_{3}^{n,\varepsilon}.
\end{align*}
Покажем, что 
\begin{equation}\label{J1to0}
\displaystyle J_{1}^{\varepsilon}=\int_{0}^{1}dx_1 \int_{0}^{1}dx_2\bigg|f^{\varepsilon}_{x_1}\bigg(x_1, x_2\bigg)-f_{x_1}\bigg(x_1,x_2\bigg)\bigg|\to 0, 
\quad \varepsilon \to 0.
\end{equation}
В самом деле, поскольку, как известно, пространство $C(K)$ всюду плотно в  $L_1(K)$, то для любого $\delta>0$ найдется такая функция $g\in C(K)$, что 
$$
\|f_{x_1} - g\|_{L_1(K)} \le \delta.
$$
В силу доказанного выше имеем также
$$
\|f^\varepsilon_{x_1} - g^\varepsilon\|_{L_1(K)} \le \delta,
$$
где $g^\varepsilon = g*\bar\varphi^\varepsilon $. Тогда находим,
\begin{align*}
&\int_{0}^{1}dx_1 \int_{0}^{1}dx_2 |f^{\varepsilon}_{x_1}(x_1, x_2) -f_{x_1}(x_1,x_2)|
 \\\\
&\le \int_{0}^{1}dx_1 \int_{0}^{1}dx_2 |f^{\varepsilon}_{x_1} (x_1, x_2) - g^\varepsilon(x_1,x_2)|
+  \int_{0}^{1}dx_1 \int_{0}^{1}dx_2 |f_{x_1} (x_1, x_2) - g(x_1,x_2)|
 \\\\
&+ \int_{0}^{1}dx_1 \int_{0}^{1}dx_2 |g^{\varepsilon} (x_1, x_2) - g(x_1,x_2)| 
 \\\\
&\le  \int_{0}^{1}dx_1 \int_{0}^{1}dx_2 |g^{\varepsilon}(x_1, x_2) - g(x_1,x_2)|  + 2\delta.
\end{align*}
Здесь в интеграле в правой части в силу равномерной непрерывности $g$ и компактности носителя функции $\varphi$ подынтегральное выражение $g^{\varepsilon}(x_1, x_2) - g(x_1,x_2)$ стремится к нулю при $\varepsilon\to 0$ при всех $(x_1,x_2)\in K$. Поскольку $\delta>0$ может быть выбрано произвольным, а интеграл $J^\varepsilon_1$ от $\delta$ не зависит, то получаем искомую сходимость (\ref{J1to0}).

\medskip

Рассмотрим второй интеграл. При фиксированном $\varepsilon$ имеем при $ n \to \infty$,
\begin{equation}\label{J2to0}
\displaystyle J_{2}^{n,\varepsilon}=\!\!\int_{0}^{1}\!dx_1 \int_{0}^{1}\!dx_2\bigg|n \left( f^{\varepsilon}\left( x_1+\frac{1}{n}, x_2\right)\!-\!f^{\varepsilon}\bigg( x_1, x_2\bigg) \right)\! -\! f^{\varepsilon}_{x_1}\bigg(x_1, x_2\bigg)\bigg| \xrightarrow{}0,
\end{equation}
по определению производной непрерывной функции $f^{\varepsilon}$ и в силу теоремы Лебега об ограниченной сходимости под знаком интеграла; здесь ограниченность конечных разностей $n \left( f^{\varepsilon}\left( x_1+\frac{1}{n}, x_2\right)\!-\!f^{\varepsilon}\bigg( x_1, x_2\bigg) \right)$ следует из теоремы о среднем, или из формулы Ньютона - Лейбница. 

\medskip

Рассмотрим $J_{3}^{n,\varepsilon}$. 
Для сокращения выкладок сперва преобразуем  подынтегральное выражение в $J_{3}^{n,\varepsilon}$:

\begin{align*} 
&n\bigg|\left( f^{\varepsilon}-f\right)\left( x_1+\frac{1}{n}, x_2\right)+ \left( f-f^{\varepsilon}\right)\bigg( x_1, x_2\bigg) \bigg|
 \\ \\
&=n\bigg|\left( f*\varphi^{\varepsilon}-f\right)\left( x_1+\frac{1}{n}, x_2\right) + \left( f-f*\varphi^{\varepsilon}\right)\bigg( x_1, x_2\bigg) \bigg|
 \\ \\
&=\bigg|\int_0^1 dz\left( f_{x_1}*\varphi^{\varepsilon}\right)\left( x_1+\frac{z}{n}, x_2\right) - \int_0^1 dz f_{x_1}\bigg( x_1+\frac{z}{n}, x_2\bigg) \bigg|.
\end{align*}
Используя теорему Фубини и формулу Ньютона -- Лейбница, сам интеграл $J_{3}^{n,\varepsilon}$ можно представить в виде 
\begin{align*}
&J_{3}^{n,\varepsilon}= \int_0^1 dx_1 \int_0^1 dx_2 \bigg|\int_{-\infty}^{\infty} dy_1 \int_{-\infty}^{\infty} dy_2 \bigg(f\bigg(x_1+\frac{1}{n}-y_1, x_2-y_2\bigg)
 \\ \\
&-f\bigg(x_1-y_1, x_2-y_2\bigg)\bigg)\displaystyle\varphi^{\varepsilon}(y_1)\varphi^{\varepsilon}(y_2)\bigg]-\bigg[f\bigg(x_1+\frac{1}{n}, x_2\bigg)- f\bigg(x_1, x_2\bigg)\bigg|
 \\ \\
&=\int_0^1 dx_1 \int_0^1 dx_2 \bigg| \int_{-\infty}^{\infty} dy_1 \int_{-\infty}^{\infty} dy_2 \bigg(f\bigg(x_1+\frac{1}{n}-y_1, x_2-y_2\bigg)\varphi^{\varepsilon}(y_1)\varphi^{\varepsilon}(y_2)
 \\ \\
&-f\bigg(x_1+\frac{1}{n}, x_2\bigg)\bigg)- \bigg( f\bigg(x_1-y_1, x_2-y_2\bigg)\varphi^{\varepsilon}(y_1)\varphi^{\varepsilon}(y_2)-f\bigg(x_1, x_2\bigg)\bigg)\bigg|. 
\end{align*}
Покажем, что это выражение стремится к нулю при $n \to \infty$ при любом фиксированном $\varepsilon>0$. Действительно, обозначим 
$$
F^\varepsilon(x_1,x_2):= \int_{-\infty}^{\infty} dy_1 \int_{-\infty}^{\infty} dy_2 f(x_1-y_1, x_2-y_2)\varphi^{\varepsilon}(y_1)\varphi^{\varepsilon}(y_2)-f(x_1, x_2).
$$
Эта функция ограничена и, стало быть, также из $L_1(K)$, и справедливо представление 
$$
J_{3}^{n,\varepsilon}=  \int_0^1 dx_1 \int_0^1 dx_2 (F^\varepsilon(x_1+\frac1{n},x_2) - F^\varepsilon(x_1,x_2)).
$$
При $n\to \infty$ последнее выражение стремится к нулю по лемме \ref{lemma0}. Стало быть, 
$$
\lim_{\varepsilon\to 0} \lim_{n\to\infty} J_{3}^{n,\varepsilon} = 0, 
$$
что и требовалось. Отсюда и из (\ref{J1to0}) и (\ref{J2to0}) утверждение (\ref{eqle1}) следует. Соотношение (\ref{le6-2}) совершенно аналогично и может быть получено транспозицией $x_1$ и $x_2$.

\medskip

\noindent
2. Докажем утверждение (\ref{eqle2}). 
Прежде всего, установим, аналогично такому же свойству для первых производных, что в условиях данного пункта леммы вторая смешанная производная $f$ также перестановочна с операцией конволюции с функцией $\varphi^\varepsilon$:
\begin{equation}\label{fxxe}
(f^\varepsilon)_{x_1x_2}(x_1,x_2) = (f_{x_1x_2}*\varphi^\varepsilon)(x_1,x_2).
\end{equation}

Далее используем тот факт, что и оператор смешанной производной второго порядка можно переставить с двойным интегралом, определяющим конволюцию, при условии абсолютной сходимости всех нужных интегралов. В самом деле, по доказанному выше в (\ref{partfe}) имеем, 
((при условии, что все интегралы абсолютно сходятся (а так оно и есть в силу условий леммы))),
\begin{align*}
&(f^\varepsilon)_{x_1x_2}(x_1,x_2) = \frac{\partial}{\partial x_2}\frac{\partial}{\partial x_1} \int_0^1 \left(\int_0^1 f(x_1-y_1,x_2-y_2)\varphi^\varepsilon(y_2) dy_2\right)\varphi^\varepsilon(y_1) dy_1
 \\\\
& =  \frac{\partial}{\partial x_2} \int_0^1 \frac{\partial}{\partial x_1}\left(\int_0^1 f(x_1-y_1,x_2-y_2)\varphi^\varepsilon(y_2) dy_2\right)\varphi^\varepsilon(y_1) dy_1.
\end{align*}
Здесь для внутреннего интеграла переменная $x_1$ является параметром, поэтому в силу абсолютной сходимости правой части верно равенство 
\begin{align*}
&  \frac{\partial}{\partial x_2} \int_0^1 \frac{\partial}{\partial x_1}\left(\int_0^1 f(x_1-y_1,x_2-y_2)\varphi^\varepsilon(y_2) dy_2\right)\varphi^\varepsilon(y_1) dy_1 
 \\\\
&= \frac{\partial}{\partial x_2} \int_0^1 \left(\int_0^1 f_{x_1}(x_1-y_1,x_2-y_2)\varphi^\varepsilon(y_2) dy_2\right)\varphi^\varepsilon(y_1) dy_1.
\end{align*}
Теперь по теореме Фубини перепишем,
\begin{align*}
&\frac{\partial}{\partial x_2} \int_0^1 \left(\int_0^1 f_{x_1}(x_1-y_1,x_2-y_2)\varphi^\varepsilon(y_2) dy_2\right)\varphi^\varepsilon(y_1) dy_1 
 \\\\
&=  \frac{\partial}{\partial x_2} \int_0^1 \left(\int_0^1 f_{x_1}(x_1-y_1,x_2-y_2)\varphi^\varepsilon(y_1) dy_1\right)\varphi^\varepsilon(y_2) dy_2,
\end{align*}
и далее, вновь используя доказанное ранее в (\ref{partfe}) и возможность дифференцировать по параметру, -- только теперь в этой роли выступит $x_2$, -- внутри интеграла, имеем,
\begin{align*}
&\frac{\partial}{\partial x_2} \int_0^1 \left(\int_0^1 f_{x_1}(x_1-y_1,x_2-y_2)\varphi^\varepsilon(y_1) dy_1\right)\varphi^\varepsilon(y_2) dy_2 
 \\\\
&= \int_0^1  \frac{\partial}{\partial x_2} \left(\int_0^1 f_{x_1}(x_1-y_1,x_2-y_2)\varphi^\varepsilon(y_1) dy_1\right)\varphi^\varepsilon(y_2) dy_2 
 \\\\
&= \int_0^1 \left(\int_0^1 f_{x_1x_2}(x_1-y_1,x_2-y_2)\varphi^\varepsilon(y_1) dy_1\right)\varphi^\varepsilon(y_2) dy_2, 
\end{align*}
что и требовалось для (\ref{fxxe}). 
Итак, 
$$
(f_{x_1x_2})^{\varepsilon}(x_1, x_2) \equiv  
(f_{x_1x_2}*\bar\varphi^\varepsilon)(x_1, x_2);
$$
будем использовать для обоих этих двух эквивалентных выражений краткое обозначение $f_{x_1x_2}^{\varepsilon}(x_1, x_2)$.

\medskip

Теперь оцениваем, 
\begin{align*}
&\int_{0}^{1}dx_1 \int_{0}^{1}dx_2 \bigg| n \left( f_{x_2}\left( x_1+\frac{1}{n}, x_2\right)
-f_{x_2}\bigg( x_1, x_2\bigg) \right) - f_{x_1x_2}\bigg(x_1, x_2\bigg)\bigg|
 \\ \\
&\leqslant\int_{0}^{1}dx_1 \int_{0}^{1}dx_2\bigg|f^{\varepsilon}_{x_1x_2}\bigg(x_1, x_2\bigg)-f_{x_1x_2}\bigg(x_1, x_2\bigg)\bigg|
 \\ \\
&+\int_{0}^{1}dx_1 \int_{0}^{1}dx_2\bigg|n \left( f^{\varepsilon}_{x_2}\left( x_1+\frac{1}{n}, x_2\right)-f^{\varepsilon}_{x_2}\bigg( x_1, x_2\bigg) \right) - f^{\varepsilon}_{x_1x_2}\bigg(x_1, x_2\bigg)\bigg|
 \\ \\
&+\int_{0}^{1}dx_1 \int_{0}^{1}dx_2\bigg|n \left( f_{x_2}\left( x_1+\frac{1}{n}, x_2\right)
-f_{x_2}\bigg( x_1, x_2\bigg) \right)
 \\ \\
&-n \left( f^{\varepsilon}_{x_2}\left( x_1+\frac{1}{n}, x_2\right)-f^{\varepsilon}_{x_2}\bigg( x_1, x_2\bigg) \right)\bigg| =J_{4}^{\varepsilon}+J_{5}^{n,\varepsilon}+J_{6}^{n,\varepsilon}
\end{align*}

Проверим сходимость к нулю каждого из интегралов в отдельности по базе $\lim_{\varepsilon\to 0} \lim_{n\to\infty}(...)$. Прежде всего, покажем, что 
\begin{equation}\label{J4to0}
\displaystyle J_{4}^{\varepsilon}=\int_{0}^{1}dx_1 \int_{0}^{1}dx_2\bigg|f^{\varepsilon}_{x_1x_2}\bigg(x_1, x_2\bigg)-f_{x_1x_2}\bigg(x_1, x_2\bigg)\bigg|\xrightarrow{}0, \quad \varepsilon \to 0.
\end{equation}
В самом деле, пусть при некотором $\delta>0$ функция $g\in C(K)$ такова, что 
$$
\|f_{x_1x_2} - g \|_{L_1(K)} \le \delta.
$$ 
Тогда в силу элементарных неравенств и доказанного выше свойства (\ref{fxxe}), а также теоремы Фубини имеем, 
$$
\|f^\varepsilon_{x_1x_2} - g^\varepsilon \|_{L_1(K)} \le \delta,
$$ 
где $g^\varepsilon = g * \varphi^\varepsilon$. В самом деле, 
\begin{align*}
&\|f^\varepsilon_{x_1x_2} - g^\varepsilon \|_{L_1(K)}  = 
\iint\limits_{K} |f^\varepsilon_{x_1x_2}(x_1,x_2) - g^\varepsilon(x_1,x_2)| dx_1dx_2
 \\\\
&=  \iint\limits_{K} |\iint\limits_{K} (f_{x_1x_2}(x_1-y_1,x_2-y_2)\! -\! g(x_1-y_1,x_2-y_2))\bar\varphi^\varepsilon(y_1,y_2)dy_1dy_2| dx_1dx_2
 \\\\
&\le  \iint\limits_{K} \left(\iint\limits_{K} |f_{x_1x_2}(x_1-y_1,x_2-y_2)\! -\! g(x_1-y_1,x_2-y_2)|\bar\varphi^\varepsilon(y_1,y_2) dx_1dx_2\right) dy_1dy_2
 \\\\
&\le  \iint\limits_{K} \left(\|f_{x_1x_2} - g \|_{L_1(K)}\right) \bar\varphi^\varepsilon(y_1,y_2) dy_1dy_2 \le \delta \iint\limits_{K}  \bar\varphi^\varepsilon(y_1,y_2)= \delta.
\end{align*}

Теперь оцениваем, при любом $\delta>0$,
\begin{align*}
&\int_{0}^{1}dx_1 \int_{0}^{1}dx_2\bigg|f^{\varepsilon}_{x_1x_2}\bigg(x_1, x_2\bigg)-f_{x_1x_2}\bigg(x_1, x_2\bigg)\bigg|
 \\\\
&\le  \int_{0}^{1}dx_1 \int_{0}^{1}dx_2\bigg|f^{\varepsilon}_{x_1x_2}\bigg(x_1, x_2\bigg)- g^\varepsilon\bigg(x_1, x_2\bigg)\bigg|
 \\\\
&+  \int_{0}^{1}dx_1 \int_{0}^{1}dx_2\bigg|f_{x_1x_2}\bigg(x_1, x_2\bigg)- g\bigg(x_1, x_2\bigg)\bigg|
 \\\\
&+ \int_{0}^{1}dx_1 \int_{0}^{1}dx_2\bigg|g^{\varepsilon}\bigg(x_1, x_2\bigg)- g\bigg(x_1, x_2\bigg)\bigg|
 \\\\
&\le   \int_{0}^{1}dx_1 \int_{0}^{1}dx_2\bigg|g^{\varepsilon}\bigg(x_1, x_2\bigg)- g\bigg(x_1, x_2\bigg)\bigg| + 2\delta.
\end{align*}
Здесь интеграл в правой части последнего неравенства стремится к нулю при $\varepsilon \to 0$, в силу равномерной непрерывности функции $g$. Поскольку же интеграл $J_{4}^{\varepsilon}$ не зависит от $\delta$, то получаем искомое стремление к нулю в (\ref{J4to0}).

\medskip

Рассмотрим член $J_{5}^{n,\varepsilon}$. При $n \to \infty$ имеем,
$$
\displaystyle J_{5}^{n,\varepsilon}=\!\!\int_{0}^{1}\!dx_1 \int_{0}^{1}\!dx_2\bigg|n \left( f^{\varepsilon}_{x_2}\left( x_1+\frac{1}{n}, x_2\right)\!-\!f^{\varepsilon}_{x_2}\bigg( x_1, x_2\bigg) \right)\! -\! f^{\varepsilon}_{x_1x_2}\bigg(x_1, x_2\bigg)\bigg| \xrightarrow{}0,
$$
согласно  определению производной для гладкой ограниченной функции $f^{\varepsilon}$ и в силу теоремы Лебега об ограниченной сходимости под знаком интеграла; здесь ограниченность конечных разностей $n \left( f^{\varepsilon}_{x_2}\left( x_1+\frac{1}{n}, x_2\right)\!-\!f^{\varepsilon}_{x_2}\bigg( x_1, x_2\bigg) \right)$ следует из теоремы о среднем, или из формулы Ньютона - Лейбница, аналогично тому, как это имело место для $J_{2}^{n,\varepsilon}$. Стало быть, получаем
$$
\displaystyle \lim_{\varepsilon\to 0} \limsup_{n\to\infty} J_{5}^{n,\varepsilon} = 0.
$$

\medskip

Рассмотрим интеграл $J_{6}^{n,\varepsilon}$. 
Преобразуем подынтегральное выражение в $J_{6}^{n,\varepsilon}$, пользуясь формулой Ньютона-Лейбница и перестановочностью операторов производной (по $x_1$) и конволюции:
\begin{align*}
&n \displaystyle \bigg| \bigg(f_{x_2}^{\varepsilon}-f_{x_2}\bigg)\left(x_1+\frac{1}{n}, x_2\right)-\bigg(f_{x_2}^{\varepsilon}-f_{x_2}\bigg)\bigg(x_1, x_2\bigg)\bigg|
 \\ \\
&=\displaystyle n  \bigg|\bigg(f_{x_2}*\varphi^{\varepsilon} - f_{x_2}\bigg) \bigg(x_1+\frac{1}{n}, x_2\bigg) + \bigg(f_{x_2}-f_{x_2}*\varphi^{\varepsilon}\bigg) \bigg(x_1, x_2\bigg)\bigg|
 \\ \\
&=\bigg|\int_0^1 dz  \bigg(f_{x_2  x_1}*\varphi^{\varepsilon}\bigg)\bigg(x_1+\frac{z}{n},x_2\bigg) - \int_0^1 dz  f_{x_2 x_1}\bigg(x_1+\frac{z}{n}, x_2\bigg)\bigg|.
\end{align*}
В результате получаем представление
\begin{align*}
&J_{6}^{n,\varepsilon}= \int_0^1 dx_1 \int_0^1 dx_2 \bigg|\int_{-\infty}^{\infty} dy_1 \int_{-\infty}^{\infty} dy_2 \bigg(f_{x_2}\bigg(x_1+\frac{1}{n}-y_1, x_2-y_2\bigg)
 \\ \\
&- f_{x_2}\bigg(x_1-y_1, x_2-y_2\bigg)\bigg)\displaystyle\varphi^{\varepsilon}(y_1)\varphi^{\varepsilon}(y_2)- \bigg(f_{x_2}\bigg(x_1+\frac{1}{n}, x_2\bigg)- f_{x_2}\bigg(x_1, x_2\bigg)\bigg)\bigg|
 \\ \\
&=\int_0^1 dx_1 \int_0^1 dx_2 \bigg| \int_{-\infty}^{\infty} dy_1 \int_{-\infty}^{\infty} dy_2 \bigg(f_{x_2}\bigg(x_1+\frac{1}{n}-y_1, x_2-y_2\bigg)\varphi^{\varepsilon}(y_1)\varphi^{\varepsilon}(y_2)
 \\ \\
&-f_{x_2}\bigg(x_1+\frac{1}{n}, x_2\bigg)\bigg)- \bigg( f_{x_2}\bigg(x_1-y_1, x_2-y_2\bigg)\varphi^{\varepsilon}(y_1)\varphi^{\varepsilon}(y_2)-f_{x_2}\bigg(x_1, x_2\bigg)\bigg)\bigg|. 
\end{align*} 
Обозначая
$$
F^\varepsilon(x_1,x_2):=\int_{-\infty}^{\infty} dy_1 \int_{-\infty}^{\infty} dy_2\bigg( f_{x_2}\bigg(x_1-y_1, x_2-y_2\bigg)\varphi^{\varepsilon}(y_1)\varphi^{\varepsilon}(y_2)-f_{x_2}\bigg(x_1, x_2\bigg)\bigg), 
$$
интеграл $J_{6}^{n,\varepsilon}$ можно переписать как 
\begin{align*}
&J_{6}^{n,\varepsilon}= \int_0^1 dx_1 \int_0^1 dx_2 \bigg| F^\varepsilon(x_1 + \frac1{n},x_2) - F^\varepsilon(x_1,x_2)\bigg|.
\end{align*} 
Поскольку функция $F^\varepsilon$ ограничена, находим, что $\lim_{n\to\infty} J_{6}^{n,\varepsilon} = 0$. Стало быть, и
$$
\lim_{\varepsilon \to 0}\lim_{n\to\infty} J_{6}^{n,\varepsilon} = 0, 
$$
что и требовалось.

\medskip

\noindent
{\bf 3.} Докажем утверждение (\ref{eqle3}): 
$$
\int_{0}^{1}dx_1 \int_{0}^{1}dx_2 \left\vert
n^2\Delta_{(x_2)}\Delta_{(x_1)}f_{}(x_1,x_2) -
n\Delta_{(x_1)}f_{x_2}(x_1,x_2)\right\vert \to 0, \quad n\to\infty. 
$$
Применяя формулу Ньютона-Лейбница в версии
$$
g(x+a) - g(x) = \int_0^1 g'(x+ sa)a ds
$$
с $a=1/n$, имеем, 

\begin{align*}
n\Delta_{(x_1)} f_{x_2}(x_1, x_2) \!=\!
n\left(f_{x_2}(x_1\!+\!\frac1{n}, x_2) \!-\! f_{x_2}(x_1, x_2)\right)
\!=\! \int_0^1 f_{x_2 x_1}(x_1\!+\!\frac{t}{n}, x_2) dt,
\end{align*}
и
\begin{align*}
\displaystyle
&n^2\Delta_{(x_2)}\Delta_{(x_1)}f(x_1, x_2) 
 \\\\
&= n^2[f(x_1+\frac1{n}, x_2+\frac1{n}) - f(x_1+\frac1{n}, x_2) - f(x_1, x_2+\frac1{n}) + f(x_1, x_2)]
 \\\\
&=  n \int_0^1 f_{x_2}(x_1+\frac1{n}, x_2+\frac{t}{n}) ds - n \int_0^1 f_{x_2}(x_1, x_2+\frac{t}{n}) dt
 \\\\
&=  n \int_0^1 [f_{x_2}(x_1+\frac1{n}, x_2+\frac{t}{n})  - f_{x_2}(x_1, x_2+\frac{t}{n})] dt
 \\\\
&=  \iint\limits_{K} f_{x_2 x_1}(x_1+\frac{s_1}{n}, x_2+\frac{s_2}{n}) ds_1 ds_2.
\end{align*}
Следовательно,
\begin{align*}
&\int_0^1 dx_1 \int_0^1 dx_2 \bigg|n\Delta_{(x_1)} f_{x_2}(x_1, x_2) - n^2\Delta_{(x_2)}\Delta_{(x_1)}f(x_1, x_2)\bigg|
 \\ \\
&= \int_0^1\! dx_1 \int_0^1\! dx_2 \bigg|\int_0^1\! f_{x_1x_2}\bigg(x_1\!+\!\frac{t_1}{n}, x_2\bigg)dt_1
\!-\! \int_0^1\! dt_1 \int_0^1\! f_{x_1x_2 }\bigg(x_1 \!+\! \frac{t_1}{n}, x_2 \!+\! \frac{t_2}{n}\bigg)dt_2\bigg|
 \\ \\
&\leqslant  \int_0^1 dx_1 \int_0^1 dx_2 \int_0^1 dt_1 \int_0^1 dt_2 \bigg|f_{x_1x_2}\bigg(x_1 +\frac{t_1}{n}, x_2\bigg) - f_{x_1x_2}\bigg(x_1+\frac{t_1}{n}, x_2+\frac{t_2}{n}\bigg)\bigg|
 \\ \\
&\leqslant  \int_0^1 dx_1 \int_0^1 dx_2 \int_0^1 dt_1 \int_0^1 dt_2 \bigg|f_{x_1x_2}\bigg(x_1+\frac{t_1}{n}, x_2 +\frac{t_2}{n}\bigg) - f_{x_1x_2}\bigg(x_1, x_2\bigg)\bigg|
 \\ \\ 
&+ \int_0^1 dx_1 \int_0^1 dx_2 \int_0^1 dt_1 \int_0^1 dt_2 \bigg|f_{x_1x_2}\bigg(x_1, x_2 \bigg) - f_{x_1x_2}\bigg(x_1+\frac{t_1}{n}, x_2\bigg)\bigg|.
\end{align*}
В силу леммы \ref{lemma0} заключаем, что оба слагаемых в правой части последнего неравенства стремятся к нулю при $n\to\infty$. Стало быть, 
$$
\int_{0}^{1}dx_1 \int_{0}^{1}dx_2 \left\vert
n^2\Delta_{(x_2)}\Delta_{(x_1)}f_{}(x_1,x_2) -
n\Delta_{(x_1)}f_{x_2}(x_1,x_2)\right\vert \to 0.
$$
Наконец, соотношение (\ref{le6-5}) эквивалентно (\ref{eqle2}), а (\ref{le6-6}) эквивалентно (\ref{eqle3}).
Лемма \ref{lemma1} доказана.
\end{proof}

\section{Доказательство теоремы \ref{thm1}}
\noindent
Ради облегчения выкладок введем обозначение:
$$
\Pi_{k_1,k_2}(x_1,x_2):=\dbinom{n-1}{k_1} \dbinom{n-1}{k_2} x_{1}^{k_1}{(1-x_1)}^{n-k_1-1} x_{2}^{k_2} {(1-x_2)}^{n-k_2-1}.
$$
Выражение в правой части естественным образом появляется при повторном дифференцировании полиномов $\tilde B_{\alpha,n}$. Отметим, что в силу формулы биному Ньютона, при всех $x_1,x_2$
\begin{equation}\label{binom}
\sum_{0\le k_1, k_2 \le n-1}\Pi_{k_1,k_2}(x_1,x_2) = 1.
\end{equation}

Покажем, что в силу соотношений (\ref{eqle2}) и (\ref{eqle3}) из леммы \ref{lemma1} следует равенство
\begin{align}\label{pr-e1}
&{\mathsf E}_{}  \int_{0}^{1} dx_1 \int_{0}^{1} dx_2 \bigg| \frac{\partial^2 \widetilde {B}_{\alpha,n}(f,x_1,x_2)}{\partial x_1 \partial x_2} - \frac{\partial^2 (f,x_1,x_2)}{\partial x_1 \partial x_2}\bigg|
  \nonumber \\  \nonumber \\  \nonumber 
&={\mathsf E}  \int_{0}^{1} dx_1 \int_{0}^{1} dx_2
\bigg| \sum\limits_{0 \leqslant k_1 < n\atop 0 \leqslant k_2 < n} n^2 \bigg[ \Delta_{(x_1)} 
\Delta_{(x_2)}f\left( \frac{k_1+\alpha_1}{n},\frac{k_2+\alpha_2}{n}\right)
 \\  \nonumber \\ 
&-\Delta_{(x_1)} \Delta_{(x_2)}f\left( x_1,x_2\right)\bigg] 
\Pi_{k_1,k_2}(x_1,x_2)\bigg|+o(1)=: I^n+o(1), \quad n\to\infty.
\end{align}
В самом деле, продифференцируем $\widetilde {B}_{\alpha,n}$ по $x_1$. Имеем, 

\begin{align*}
&\frac{\partial \widetilde {B}_{\alpha,n}(f,x_1,x_2)}{\partial x_1}=
\frac{\partial}{\partial x_1} \sum\limits_{k_1, k_2=0}^{n} f\left( \frac{k_1}{n}+\frac{\alpha_1}{n},\frac{k_2}{n}+\frac{\alpha_2}{n}\right) \dbinom{n}{k_1} \dbinom{n}{k_2} 
 \\ \\
&\times x_{1}^{k_1} {(1-x_1)}^{n-k_1} x_{2}^{k_2} {(1-x_2)}^{n-k_2}
= \sum\limits_{0 < k_1 \leqslant n\atop 0 \leqslant k_2 \leqslant n}
f\left( \frac{k_1}{n}+\frac{\alpha_1}{n},\frac{k_2}{n}+\frac{\alpha_2}{n}\right)  \dbinom{n}{k_1} \dbinom{n}{k_2}
 \\ \\
&\times k_{1} x_{1}^{k_1-1} {(1-x_1)}^{n-k_1} x_{2}^{k_2} {(1-x_2)}^{n-k_2}
- \sum\limits_{0 \leqslant k_1 < n\atop 0 \leqslant k_2 \leqslant n}
f\left( \frac{k_1}{n}+\frac{\alpha_1}{n},\frac{k_2}{n}+\frac{\alpha_2}{n}\right)  \dbinom{n}{k_1} \dbinom{n}{k_2}
 \\ \\
&\times x_{1}^{k_1} (n-k_1) {(1-x_1)}^{n-k_1-1} x_{2}^{k_2} {(1-x_2)}^{n-k_2}
= \sum\limits_{0 \leqslant k_1 < n\atop 0 \leqslant k_2 \leqslant n}
f\left( \frac{k_1+\alpha_1 + 1}{n},\frac{k_2+\alpha_2}{n}\right) 
 \\ \\
&\times \dbinom{n}{k_1 + 1} \dbinom{n}{k_2} (k_{1} + 1) x_{1}^{k_1}
{(1-x_1)}^{n-k_1-1} x_{2}^{k_2} {(1-x_2)}^{n-k_2}
 \\ \\
&- \sum\limits_{0 \leqslant k_1 < n\atop 0 \leqslant k_2 \leqslant n}
f\left( \frac{k_1+\alpha_1}{n},\frac{k_2+\alpha_2}{n}\right) 
\dbinom{n}{k_1} \dbinom{n}{k_2} x_{1}^{k_1} (n-k_1)
{(1-x_1)}^{n-k_1-1} x_{2}^{k_2} {(1-x_2)}^{n-k_2}
 \\ \\
&=\sum\limits_{0 \leqslant k_1 < n\atop 0 \leqslant k_2 \leqslant n}
n \Delta_{(x_1)} f\left( \frac{k_1+\alpha_1}{n},\frac{k_2+\alpha_2}{n}\right) 
\dbinom{n-1}{k_1} \dbinom{n}{k_2} x_{1}^{k_1}
{(1-x_1)}^{n-k_1-1} x_{2}^{k_2} {(1-x_2)}^{n-k_2}.
\end{align*}
Продифференцируем теперь полученное выражение по $x_2$:
\begin{align*}
&\frac{\partial^2 \widetilde {B}_{\alpha,n}(f,x_1,x_2)}{\partial x_1 \partial x_2}=
\frac{\partial^2}{\partial x_1 \partial x_2} \sum\limits_{k_1, k_2=0}^{n} f\left( \frac{k_1+\alpha_1}{n},\frac{k_2+\alpha_2}{n}\right) \dbinom{n}{k_1} \dbinom{n}{k_2}&&&\\ \\
&\times x_{1}^{k_1} {(1-x_1)}^{n-k_1} x_{2}^{k_2} {(1-x_2)}^{n-k_2}
=\frac{\partial}{\partial x_2} \sum\limits_{0 \leqslant k_1 < n\atop 0 \leqslant k_2 \leqslant n} n \Delta_{(x_1)} f\left( \frac{k_1+\alpha_1}{n},\frac{k_2+\alpha_2}{n}\right)&&&\\ \\ 
&\times \dbinom{n-1}{k_1} \dbinom{n}{k_2} \times x_{1}^{k_1} {(1-x_1)}^{n-k_1-1} x_{2}^{k_2} {(1-x_2)}^{n-k_2}&&&\\ \\
&= \sum\limits_{0 < k_1 \leqslant n\atop 0 \leqslant k_2 < n}
n \Delta_{(x_1)} f\left( \frac{k_1+\alpha_1}{n},\frac{k_2+\alpha_2}{n}\right)  \dbinom{n-1}{k_1} \dbinom{n}{k_2}
x_{1}^{k_1} {(1-x_1)}^{n-k_1-1} k_2 x_{2}^{k_2-1} {(1-x_2)}^{n-k_2}&&&
 \\ \\
&-\! \sum\limits_{0 \leqslant k_1 < n\atop 0 < k_2 \leqslant n}\! n \Delta_{(x_1)}
f\left( \frac{k_1+\alpha_1}{n},\frac{k_2+\alpha_2}{n}\right)  \dbinom{n-1}{k_1} \dbinom{n}{k_2}&&&
 \\ \\
&\times x_{1}^{k_1} {(1-x_1)}^{n-k_1-1} x_{2}^{k_2}(n-k_2) {(1-x_2)}^{n-k_2-1}
= \sum\limits_{0 \leqslant k_1 < n\atop 0 \leqslant k_2 < n}
f\left( \frac{k_1+\alpha_1 }{n},\frac{k_2+\alpha_2+ 1}{n}\right) &&&
 \\ \\
&\times \dbinom{n-1}{k_1} \dbinom{n}{k_2+1} x_{1}^{k_1}
{(1-x_1)}^{n-k_1-1} x_{2}^{k_2} {(1-x_2)}^{n-k_2-1}&&&
 \\ \\
&- \!\!\!\!\sum\limits_{0 \leqslant k_1 < n\atop 0 < k_2 \leqslant n}\!\!\!n \Delta_{(x_1)}
f\left( \frac{k_1+\alpha_1}{n},\frac{k_2+\alpha_2}{n}\right)  \dbinom{n-1}{k_1} \dbinom{n}{k_2} 
x_{1}^{k_1} {(1-x_1)}^{n-k_1-1} x_{2}^{k_2}(n-k_2) {(1-x_2)}^{n-k_2-1}&&&
 \\ \\
&=\sum\limits_{0 \leqslant k_1 < n\atop 0 \leqslant k_2 < n}
n^2 \left( \Delta_{(x_1)} f\left( \frac{k_1+\alpha_1}{n},\frac{k_2+\alpha_2 +1}{n}\right)
-\Delta_{(x_1)} f\left( \frac{k_1+\alpha_1}{n},\frac{k_2+\alpha_2}{n}\right) \right)&&&
 \\ \\
&\times \dbinom{n-1}{k_1} \dbinom{n-1}{k_2} x_{1}^{k_1}
{(1-x_1)}^{n-k_1-1} x_{2}^{k_2} {(1-x_2)}^{n-k_2-1}&&&
 \\\\
&=\sum\limits_{0 \leqslant k_1 < n\atop 0 \leqslant k_2 < n}
n^2 \Delta_{(x_1)} \Delta_{(x_2)}f \left( \frac{k_1+\alpha_1}{n},\frac{k_2+\alpha_2}{n}\right)
\dbinom{n-1}{k_1} \dbinom{n-1}{k_2}
 \\ \\
&\times x_{1}^{k_1}
{(1-x_1)}^{n-k_1-1} x_{2}^{k_2} {(1-x_2)}^{n-k_2-1}.
\end{align*}

Теперь  формула (\ref{pr-e1}) следует из найденного представления для смешанной производной $\frac{\partial^2 \widetilde {B}_{\alpha,n}(f,x_1,x_2)}{\partial x_1 \partial x_2}$  и из леммы \ref{lemma1}.

\medskip

Оценим сверху (неотрицательный) интеграл $I^n$ в правой части(\ref{pr-e1}), используя сглаженную функцию $f^\varepsilon = f*\bar \varphi^\varepsilon$, ранее введенную в доказательстве леммы \ref{lemma1}. Имеем,

\begin{eqnarray*}
\label{first}
&I^n\leqslant n^2 {\mathsf E}_{}  \displaystyle \int_{0}^{1} dx_1 \int_{0}^{1} dx_2     \bigg|\sum\limits_{0 \leqslant k_1 < n\atop 0 \leqslant k_2 < n}
\Biggm(\Delta_{(x_1)}\Delta_{(x_2)}f^{\varepsilon}\left( \frac{k_1+\alpha_1}{n},\frac{k_2+\alpha_2}{n}\right)\nonumber \\ \nonumber\\
&-\Delta_{(x_1)}\Delta_{(x_2)} f^{\varepsilon}\left( x_1,x_2 \right) \Biggm)\Pi_{k_1,k_2}(x_1,x_2)\bigg|
 \\\\
\label{second}
&+n^2 {\mathbb E}_{}  \displaystyle \int_{0}^{1} dx_1 \int_{0}^{1} dx_2 
\bigg|\sum\limits_{0 \leqslant k_1 < n\atop 0 \leqslant k_2 < n}
\Biggm(\Delta_{(x_1)}\Delta_{(x_2)}f^{\varepsilon}\left(x_1,x_2\right)-\Delta_{(x_1)} \Delta_{(x_2)}f\left( x_1,x_2 \right)\Biggm) \nonumber \\ \nonumber\\
&\times\Pi_{k_1,k_2}(x_1,x_2)\bigg|
 \\\\
\label{third}
&+n^2 {\mathbb E}_{}  \displaystyle \int_{0}^{1} dx_1 \int_{0}^{1} dx_2 
\bigg|\sum\limits_{0 \leqslant k_1 < n\atop 0 \leqslant k_2 < n}
\Biggm( \Delta_{(x_1)}\Delta_{(x_2)}f^{\varepsilon}\left( \frac{k_1+\alpha_1}{n},\frac{k_2+\alpha_2}{n}\right)\nonumber \\ \nonumber\\
&-\displaystyle \Delta_{(x_1)}\Delta_{(x_2)} f\left( \frac{k_1+\alpha_1}{n},\frac{k_2+\alpha_2}{n} \right)\Biggm)\Pi_{k_1,k_2}(x_1,x_2)\bigg|
  \\\\
&=: I_1^{n,\varepsilon}+I_2^{n,\varepsilon}+I_3^{n,\varepsilon}.    
\end{eqnarray*}
Все три выражения $I_k^{n,\varepsilon}, 1\le k\le 3$, неотрицательны.  
Покажем, что 
$$
\lim_{\varepsilon\to 0} 
\limsup_{n\to\infty} I_k^{n,\varepsilon} = 0, \qquad 
k=1,2,3.
$$
Здесь удобно $I_2^{n,\varepsilon}$ рассмотреть в последнюю очередь, см. конец доказательства.

Рассмотрим интеграл $I_2^{n,\varepsilon}$. 
Оценим его сверху, разложив, в свою очередь, оценку в сумму трех интегралов: в силу (\ref{binom}) имеем, 
\begin{align*}
&I_2^{n,\varepsilon} 
= \int_0^1 dx_1 \int_0^1 dx_2 \bigg|n^2\Delta_{(x_1)}\Delta_{(x_2)}f^{\varepsilon}\left(x_1, x_2\right)
-n^2\Delta_{(x_1)}\Delta_{(x_2)}f\left(x_1, x_2\right)\bigg|
 \nonumber \\\nonumber \\
& 
\leqslant \int_0^1 dx_1 \int_0^1 dx_2 \bigg|n^2 \Delta_{(x_1)}\Delta_{(x_2)}f^{\varepsilon}(x_1, x_2)-f_{x_1 x_2}^{\varepsilon}(x_1, x_2)\bigg|
 \\\nonumber \\  
&+\int_0^1 dx_1 \int_0^1 dx_2 \bigg| f_{x_1 x_2}^{\varepsilon}(x_1, x_2) -f_{x_1 x_2}(x_1, x_2)\bigg|
 \\\nonumber\\  
&+\int_0^1 dx_1 \int_0^1 dx_2 \bigg|f_{x_1 x_2}(x_1, x_2)-n^2\Delta_{(x_1)}\Delta_{(x_2)}f(x_1, x_2)\bigg|
 \\\nonumber \\\nonumber 
&=:I_{21}^{n,\varepsilon}+I_{22}^{\varepsilon}+I_{23}^{n}.
\end{align*}
Рассмотрим интеграл $I_{22}^{\varepsilon}$. Поскольку $f_{x_1x_2} \in L_1(K)$, то он сходится к нулю при $\varepsilon \to 0$, в силу леммы \ref{lemma0}. 

Рассмотрим интеграл $I_{21}^{n,\varepsilon}$. Оценим его таким образом: 
\begin{align*}
0\le  I_{21}^{n,\varepsilon} &\le
 \int_0^1 dx_1 \int_0^1 dx_2 \bigg|n^2 \Delta_{(x_1)}\Delta_{(x_2)}f^{\varepsilon}(x_1, x_2) - n \Delta_{(x _1)}f_{x_2}^{\varepsilon}(x_1, x_2)\bigg|
  \\\\
&+ \int_0^1 dx_1 \int_0^1 dx_2 \bigg|n \Delta_{(x _1)}f_{x_2}^{\varepsilon}(x_1, x_2)-f_{x_1 x_2}^{\varepsilon}(x_1, x_2)\bigg|  =: I_{21a}^{n,\varepsilon}+ I_{21b}^{n,\varepsilon}.
\end{align*}
При $n\to\infty$ имеем сходимость $I_{{21}b}^{n,\varepsilon} \to 0$ 
в силу формулы (\ref{eqle2}) леммы \ref{lemma1}. 
Аналогично, имеем
$\lim_{n\to\infty} I_{21a}^{n,\varepsilon} = 0$ в силу формулы (\ref{le6-6}) леммы \ref{lemma1}.

\medskip

Рассмотрим интеграл $I_{23}^{n}$. Оценим его сверху суммой двух интегралов. Поскольку 
$\Delta_{(x_2)}\Delta_{(x_1)}f(x_1, x_2) = \Delta_{(x_1)}\Delta_{(x_2)}f(x_1, x_2)$,
то  оценим,
\begin{align*}
&\int_0^1 dx_1 \int_0^1 dx_2 \bigg|f_{x_1 x_2}(x_1, x_2)-n^2\Delta_{(x_1)}\Delta_{(x_2)}f(x_1, x_2)\bigg|\leqslant
 \\ \\
&\leqslant\int_0^1 dx_1 \int_0^1 dx_2 \bigg|n \Delta_{(x_1)} f_{x_2}(x_1, x_2)-n^2\Delta_{(x_2)}\Delta_{(x_1)}f(x_1, x_2)\bigg|
 \\ \\
&+\int_0^1 dx_1 \int_0^1 dx_2 \bigg|f_{x_1 x_2}(x_1, x_2)-n \Delta_{(x_1)} f_{x_2}(x_1, x_2)\bigg|= I_{23a}^{n}+I_{23b}^{n}.
\end{align*}

Для  $I_{23a}^{n}$ имеем, $I_{23a}^{n} \xrightarrow{}0$ при $n \to \infty$ в силу соотношения (\ref{eqle3}) леммы \ref{lemma1}.

\medskip

Для  $I_{23b}^{n}$ имеем $I_{23b}^{n} \xrightarrow{}0$ при $n \to \infty$ в силу соотношения (\ref{eqle2}) леммы \ref{lemma1}. 

\medskip

Таким образом, для интеграла $I_2^{n, \varepsilon}$  получаем 
$$
0\leqslant \lim_{\varepsilon\to 0}
\limsup_{n\to\infty} 
I_2^{n,\varepsilon}\leqslant \lim_{\varepsilon\to 0}
\limsup_{n\to\infty} (I_{21}^{n, \varepsilon}+I_{22}^{\varepsilon}+I_{23}^{n}) = 0.
$$

Рассмотрим теперь интеграл $I_3^{n, \varepsilon}$. 
Напомним определение и свойство бэта-функции, которое будет полезно в дальнейшем:
$$
\displaystyle \mathcal{B}(x, y)  = \int_{0}^{1} t^{x-1}\left( 1-t\right) ^{y-1} dt =\frac{\Gamma(x)\Gamma(y)}{\Gamma(x+y)}.
$$
\newline
Рассмотрим каждый из членов двойной суммы в этом интеграле отдельно без множителя $\Pi_{k_1,k_2}(x_1,x_2)$, причем внесем модуль внутрь суммы, отчего неравенство только усилится. 
При любом $0\le k_1,k_2<n$ имеют место равенства
\begin{align*}
&\int_{0}^{1} d{\alpha_1} \int_{0}^{1} d{\alpha_2}
\bigg|\Delta_{(x_1)}\Delta_{(x_2)}f^{\varepsilon}\left( \frac{\alpha_1+k_1}{n},\frac{\alpha_2+k_2}{n}\right)-
\Delta_{(x_1)} \Delta_{(x_2)}f\left( \frac{\alpha_1+k_1}{n},\frac{\alpha_2+k_2}{n} \right) \bigg|
 \\ \\
&=n^2 \int_{k_1/n}^{(k_1+1)/n} d{a_{1}} \int_{k_2/n}^{(k_2+1)/n} d{a_{2}}
\bigg|\Delta_{(x_1)}\Delta_{(x_2)}f^{\varepsilon}\left( a_{1},a_{2}\right)-\Delta_{(x_1)}\Delta_{(x_2)} f\left( \alpha_{1}^{'},\alpha_{2}^{'} \right) \bigg|.
\end{align*}
Складывая по всем $0\le k_1,k_2<n$, в итоге получаем, 
\begin{align}\label{*}
&\sum\limits_{0 \leqslant k_1 < n\atop 0 \leqslant k_2 < n}
\int_{0}^{1} d{\alpha_1} \int_{0}^{1} \bigg|
\Delta_{(x_1)}\Delta_{(x_2)}f^{\varepsilon}\left( \frac{k_1+\alpha_1}{n},\frac{k_2+\alpha_2}{n}\right)
 \nonumber \\ \nonumber \\
&-\Delta_{(x_1)}\Delta_{(x_2)} f\left( \frac{k_1+\alpha_1}{n},\frac{k_2+\alpha_2}{n} \right) d{\alpha_2} \bigg|
 \\ \nonumber \\ \nonumber
&=n^2 \int_{0}^{\frac{n-1}{n}}d\alpha_{1}^{'}
\int_{0}^{\frac{n-1}{n}}da_{2}\bigg| \Delta_{(x_1)}\Delta_{(x_2)}f^{\varepsilon}\left( a_{1},a_{2}\right)-\Delta_{(x_1)}\Delta_{(x_2)} f\left( a_{1},a_{2}\right) \bigg|. 
\end{align}
Данное тождество будет использовано в конце выкладок.

\medskip

Теперь приступим к преобразованию интеграла $I_3^{n, \varepsilon}$, используя известные комбинаторные выражения для гамма-функции от целочисленного аргумента через факториалы. Как вскоре выяснится, с их помощью для $I_3^{n, \varepsilon}$ получается представление, почти в точности эквивалентное интегралу $I_2^{n, \varepsilon}$. Имеем,
\begin{align*}
&I_3^{n, \varepsilon}=n^2 {\mathsf E}_{} \int_0^1 dx_1 \int_0^1 \bigg|\bigg[\sum\limits_{0 \leqslant k_1 < n\atop 0 \leqslant k_2 < n}
\Delta_{(x_1)}\Delta_{(x_2)}f^{\varepsilon}\left(\frac{k_1+\alpha_1}{n},
\frac{k_2+\alpha_2}{n}\right)
 \\ \\
&-\Delta_{(x_1)}\Delta_{(x_2)}f\left(\frac{k_1+\alpha_1}{n}, \frac{k_2+\alpha_2}{n}\right) \bigg]
\end{align*}
\begin{align*}
&\times x^{k_1}(1-x)^{n-k_1-1}x^{k_2}(1-x)^{n-k_2-1}\left(\begin{array}{c}n-1 \\ k_{1}\end{array}\right)\left(\begin{array}{c}n-1 \\ k_{2}\end{array}\right) \bigg|dx_2
 \\ \\
&\leqslant n^2 {\mathsf E}_{} \bigg[\sum\limits_{0 \leqslant k_1 < n\atop 0 \leqslant k_2 < n}
\int_0^1 dx_1 \int_0^1 \bigg|
\Delta_{(x_1)}\Delta_{(x_2)}f^{\varepsilon}\left(\frac{k_1+\alpha_1}{n},\frac{k_2+\alpha_2}{n}\right)
 \\ \\
&-\Delta_{(x_1)}\Delta_{(x_2)}f\left(\frac{k_1+\alpha_1}{n}, \frac{k_2+\alpha_2}{n}\right)\bigg]
 \\ \\
&\times x^{k_1}(1-x)^{n-k_1-1}x^{k_2}(1-x)^{n-k_2-1}\left(\begin{array}{c}n-1 \\ k_{1}\end{array}\right)\left(\begin{array}{c}n-1 \\ k_{2}\end{array}\right) dx_2 \bigg|
\end{align*}
\begin{align*}
&=n^2 {\mathsf E}_{} \sum\limits_{0 \leqslant k_1 <n\atop 0 \leqslant k_2< n} \bigg|
\Delta_{(x_1)}\Delta_{(x_2)}f^{\varepsilon}\left(\frac{k_1+\alpha_1}{n},\frac{k_2+\alpha_2}{n}\right)
-\Delta_{(x_1)}\Delta_{(x_2)}f\left(\frac{k_1+\alpha_1}{n}, \frac{k_2+\alpha_2}{n}\right)\bigg| \times
 \\ \\
&\times\int_0^1 dx_1 \int_0^1 \bigg[ x^{k_1}(1-x)^{n-k_1-1}x^{k_2}(1-x)^{n-k_2-1}\left(\begin{array}{c}n-1 \\ k_{1}\end{array}\right)\left(\begin{array}{c}n-1 \\ k_{2}\end{array}\right) dx_2 \bigg]
\end{align*}
\begin{align*}
&= {\mathsf E}_{} \sum\limits_{0 \leqslant k_1 < n\atop 0 \leqslant k_2 < n}
\bigg| \Delta_{(x_1)}\Delta_{(x_2)}f^{\varepsilon}\left(\frac{k_1+\alpha_1}{n},\frac{k_2+\alpha_2}{n}\right)-
\Delta_{(x_1)}\Delta_{(x_2)}f\left(\frac{k_1+\alpha_1}{n}, \frac{k_2+\alpha_2}{n}\right) \bigg|
 \\ \\
&\times n^2 \left(\begin{array}{c}n-1 \\ k_{1}\end{array}\right)\left(\begin{array}{c}n-1 \\ k_{2}\end{array}\right)\mathcal{B}(k_1+1, n-k_1)\mathcal{B}(k_2+1, n-k_2)
 \\ \\
&= {\mathsf E}_{} \sum\limits_{0 \leqslant k_1 < n\atop 0 \leqslant k_2 < n}
\bigg| \Delta_{(x_1)}\Delta_{(x_2)}f^{\varepsilon}\left(\frac{k_1+\alpha_1}{n},\frac{k_2+\alpha_2}{n}\right)-
\Delta_{(x_1)}\Delta_{(x_2)}f\left(\frac{k_1+\alpha_1}{n}, \frac{k_2+\alpha_2}{n}\right) \bigg|
 \\ \\
&\times n^2\left(\begin{array}{c}n-1 \\ k_{1}\end{array}\right)\left(\begin{array}{c}n-1 \\ k_{2}\end{array}\right)\frac{\Gamma(k_1+1)\Gamma(n-k_1)}{\Gamma(n+1)}
\frac{\Gamma(k_2+1)\Gamma(n-k_2)}{\Gamma(n+1)}
\end{align*}
\begin{align*}
&={\mathsf E}_{} \sum\limits_{0 \leqslant k_1 < n\atop 0 \leqslant k_2 < n}
\bigg| \Delta_{(x_1)}\Delta_{(x_2)}f^{\varepsilon}\left(\frac{k_1+\alpha_1}{n},\frac{k_2+\alpha_2}{n}\right)-
\Delta_{(x_1)}\Delta_{(x_2)}f\left(\frac{k_1+\alpha_1}{n}, \frac{k_2+\alpha_2}{n}\right) \bigg|
 \\ \\
&\times\left(\begin{array}{c}n-1 \\ k_{1}\end{array}\right)\left(\begin{array}{c}n-1 \\ k_{2}\end{array}\right)\frac{k_1!(n-k_1-1)!}{n!} \frac{k_2!(n-k_2-1)!}{(n)!}
 \\ \\
&={\mathsf E}_{} \sum\limits_{0 \leqslant k_1 < n\atop 0 \leqslant k_2 < n}
\bigg| \Delta_{(x_1)}\Delta_{(x_2)}f^{\varepsilon}\left(\frac{k_1+\alpha_1}{n},\frac{k_2+\alpha_2}{n}\right)-
\Delta_{(x_1)}\Delta_{(x_2)}f\left(\frac{k_1+\alpha_1}{n}, \frac{k_2+\alpha_2}{n}\right) \bigg|
 \\ \\
&\stackrel{(\ref{*})}{=}n^2 \int_{0}^{1-\frac{1}{n}} da_{1} 
\int_{0}^{1-\frac{1}{n}} da_{2}
\bigg|\Delta_{(x_1)}\Delta_{(x_2)}f^{\varepsilon}\left(a_{1},a_{2}\right)-
\Delta_{(x_1)}\Delta_{(x_2)}f\left(a_{1}, a_{2}\right)\bigg|.
\end{align*}
Из вида последнего выражения ясно видно, что сходимость интеграла 
$$
\lim_{\varepsilon\to 0} \limsup_{n\to\infty}I_3^{n, \varepsilon} = 0
$$ 
следует буквально из тех же выкладок, что и для $I_2^{n, \varepsilon}$, с заменой $x_i$ на $a_i$, $i=1,2$.

\medskip

Наконец, сходимость $\lim_{n\to\infty}I_1^{n, \varepsilon} = 0$ {\bf при выбранном ранее для интегралов $I_2^{n, \varepsilon}$ и $I_3^{n, \varepsilon}$ сколь угодно малом $\varepsilon>0$}  вытекает из предложения \ref{Veretennikov} при $d=2$ и для мультииндекса $k=(1,1)$, в силу гладкости $f^\varepsilon$. 

\medskip

Таким образом, 
$$
\lim_{\varepsilon\to 0} \limsup_{n\to\infty}(
I_1^{n, \varepsilon} 
+ I_2^{n, \varepsilon}
+ I_3^{n, \varepsilon}) = 0.
$$ 
Стало быть, 
$$
\lim_{n\to\infty}I^{n} =0, 
$$
что и доказывает искомую сходимость
$$
\displaystyle {\mathsf E}_{} \int\limits_0^1 dx_1 \int\limits_0^1 dx_2 \bigg| \frac{\partial^2 \widetilde{B}_{\alpha,n}(f,x_1,x_2)}{\partial x_1 \partial x_2} - \frac{\partial^2 f(x_1,x_2)}{\partial x_1 \partial x_2}\bigg| \xrightarrow{} 0, \quad n \to \infty.$$
Теорема \ref{thm1} доказана. \hfill $\square$

\section{Замечания о смешанных производных}
Приведем несколько теорем  о  смешанных  производных, формально отличных, но близких по духу к теореме \ref{thm1}, а также имеющий к ним отношение один результат о первой производной В.В. Степанова.

\begin{Proposition}[Schwarz, цит. по  {\cite[теорема 8.2.3]{Zorich}}]\label{Schwarz}
Если функция $f: G \mapsto \mathbb R$ имеет в области $G$ частные производные $f_{xy}$ и $f_{yx}$, то в любой точке $(x_0,y_0)\in G$, в которой обе они непрерывна, они равны:
$$
f_{xy}(x_0,y_0) = f_{yx}(x_0,y_0).
$$
\end{Proposition}
\noindent
Полезное обобщение см. в  предложении \ref{zorich-ex}, взятом из упражнения 8.4.2b там же. В нем не предполагается заранее существование обеих частных производных в точке, а лишь одной из них; существование второй следует из утверждения.
В следующей также классической теореме предполагается существование второго дифференциала в точке. 

\begin{Proposition}[Young {\cite{Young}}]\label{young}
Рассмотрим функцию $\displaystyle f: D \subset \mathbf{R}^{2} \rightarrow \mathbf{R}$ и некоторую точку $(x^0_1, x^0_2) \in D$. Если частные производные функции первого порядка $\displaystyle \partial f/\partial x_1$ и $\displaystyle \partial f/\partial x_2$ 
существуют
в некоторой окрестности точки $(x^0_1, x^0_2)$ и дифференцируемы в точке $(x^0_1, x^0_2)$, то
значения обеих смешанных частных производных второго порядка функции $f$ в точке $(x^0_1, x^0_2)$ равны: 
$$
\displaystyle\left(\frac{\partial^2f}{\partial x_1  \partial x_2}\right)_{\left(x^0_1, x^0_2\right)}=\left(\frac{\partial^2f}{\partial x_2 \partial x_1}\right)_{\left(x^0_1, x^0_2\right)}.
$$
\end{Proposition}
В оригинальной статье Янга 1908 года многостраничное доказательство читаемо поистине с трудом, как порой случается с первыми доказательствами. Современное изящное и короткое обоснование можно найти в {\cite[(8.12.3)]{Dieudonne}, а в интернете одностраничная версия Дьёдонне доступна еще по линку \cite{Morrow}. 
В учебных целях процитируем обе формулировки,  сохраняя язык  оригиналов.

\begin{Proposition}[Young -- Dieudonne]
(A)
{\cite[(8.12.3)]{Dieudonne}}
Let $G$ be an open set in $\mathbb R^n$; if a mapping $f$ of $G$ into a Banach space $F$ is twice differentiable at $x_0$, then the partial derivatives $D_if$ are differentiable at $x_0$ and 
$$
D_iD_jf(x_0) = D_jD_if(x_0).
$$

\noindent
(B)\cite{Morrow} 
Suppose $f(x, y)$ is defined in a neighborhood of a point $(a, b)$. Suppose the partial derivatives
$f_x, f_y$ are defined in a neighborhood of $(a, b)$ and are differentiable at $(a, b)$. (In particular this implies
that $f_x, f_y$ are continuous at $(a, b)$, but it is not assumed that their derivatives exist anywhere other than
at $(a, b)$.) A short statement of the assumption is that $Df = [f_x, f_y]$ is differentiable at $(a, b)$. This is
sometimes stated as f is twice differentiable at $(a, b)$. Then
$$
(f_x)_y(a, b) = (f_y)_x(a, b),
$$
sometimes stated as
$$
f_{xy}(a, b) = f_{yx}(a, b).
$$

\end{Proposition}

\begin{Proposition}[Aksoy, Martelli \cite{Aksoy}]\label{Fubini}
Следующие утверждения эквивалентны.\\ \\
(1) \quad Путь $g \in C(U, \mathbb R^2)$ и $[a, b] * [c, d] \subset U$. Тогда
$$
\int_a^b \int_c^d g(x, y)dydx = \int_c^d \int_a^b g(x, y)dxdy.
$$
(2) \quad Пусть $f \in C^1(U, \mathbb R^2)$ и существует смешанная производная $f_{xy} \in C(U, R^2)$. Тогда производная $f_{yx}$ также существует и 
$$
f_{yx}(x,y) = f_{xy} (x,y) \quad \text{при всех  \; $(x,y) \in U$}.
$$
\end{Proposition}

\begin{Proposition}[Stepanoff \cite{Stepanoff}, {\cite[Теорема 3.1.9]{Federer}}]
Если $\displaystyle E$ -- измеримое множество,  $E\subset G \subset \mathbf{R}^{m}$, где  множество $G$ открыто, и 
$\displaystyle f: G \rightarrow \mathbf{R}^{n}$ локально ограничена, измерима, и $\displaystyle\limsup\limits_{x\to a}\frac{|f(x)-f(a)|}{|x-a|} < \infty$ для почти всех  точек $a\!\in\! E$, то функция $f$ дифференцируема в $L_m$-- почти всех точках множества $E$, где $L_m$- m-мерная мера Лебега.
\end{Proposition}
Конечно, эта последняя теорема применима и ко вторым производным. Правда, вопроса о равенстве двух смешанных производных она не касается.

\begin{Proposition}[{\cite[упражнение 8.4.2b]{Zorich}}]\label{zorich-ex}
Пусть функция $f$ имеет частные производные $f_x$, $f_y$ в некоторой окрестности точки $(x_0,y_0)$. Тогда, если смешанная производная $f_{xy}$ (или $f_{yx}$) существует в $U$ и непрерывна в $(x_0,y_0)$, то смешанная производная $f_{yx}$ (соответственно, $f_{xy}$) также существует в этой точке и имеет место равенство 
$$
f_{yx}(x_0,y_0) = f_{xy} (x_0,y_0).
$$
\end{Proposition}


\begin{Proposition}[Толстов {\cite[Теорема 7]{Tolstov}}]
Если $f(x_1, x_2)$ линейно непрерывна (то есть, непрерывна по каждому аргументу), измерима в области G и в точках множества Е плоской положительной меры для каждой из ее частных производных $\displaystyle \frac{\partial f(x_1,x_2)}{\partial x_1}$ и 
$\displaystyle \frac{\partial f(x_1,x_2)}{\partial x_2}$ все производные числа по каждому из переменных конечны, то почти всюду на Е существуют и совпадают между собой смешанные производные $\displaystyle \frac{\partial^2 f(x_1,x_2)}{\partial x_1 \partial x_2}$ и $\displaystyle \frac{\partial^2 f(x_1,x_2)}{\partial x_2 \partial x_1}$.
\end{Proposition}

\begin{Proposition}[Толстов {\cite[Теорема 8]{Tolstov}}]\label{thm2}
Если функция $f(x_1, x_2)$ линейно непрерывна, всюду в области $G$ допускает производные,  $\displaystyle \frac{\partial^2 f(x_1, x_2)}{\partial x^2_1}, \quad$ 
$\displaystyle \frac{\partial^2 f(x_1, x_2)}{\partial x_1 \partial x_2}, \quad$
$\displaystyle \frac{\partial^2 f(x_1, x_2)}{\partial x_2 \partial x_1}, \quad$
$\displaystyle \frac{\partial^2 f(x_1, x_2)}{\partial x^2_2},$ 
то почти всюду в G:
$$
\displaystyle \frac{\partial^2 f(x_1, x_2)}{\partial x_1 \partial x_2} = \displaystyle \frac{\partial^2 f(x_1, x_2)}{\partial x_2 \partial x_1}.
$$
    
\end{Proposition}
 
Наконец, совсем просто (по модулю теорем Шварца, или Янга -- см. предложения \ref{Schwarz} и \ref{young}) вопрос  о равенстве двух смешанных производных решается в случае соболевских производных. Для удобства читателя напомним, что функция $g$ является смешанной соболевской производной второго порядка по переменым $x$ и $y$ в $L_p(G)$, если найдется последовательность гладких функций $f^n \in C^\infty(G)$  таких, что 
$$
\|f^n - f\|_{L_p(G)} \to 0, \; n\to\infty, 
$$
и при этом 
$$
\|f^n_{xy} - g\|_{L_p(G)} \to 0, \; n\to\infty.
$$
Следующий элементарный результат вытекает из того факта, что для гладких функций  $f^n_{xy} = f^n_{yx}$ во всех точках. 
\begin{Proposition}[О cоболевских производных]\label{Sobolev}   
Если функция $f$ из класса $L_p(G)$ с 
$p>0$ обладает смешанной соболевской производной $f_{xy}\in L_p(G)$, то также существует  соболевская производная $f_{yx}\in L_p(G)$, и функции $f_{yx}$ и  $f_{xy}$ почти всюду в $G$  равны: 
$
\|f_{yx} - f_{xy}\|_{L_p(G)} = 0
$.
\end{Proposition}

\end{document}